\documentclass{scrartcl}

\usepackage{graphicx}        
\usepackage{lmodern}

\usepackage{amssymb,amsmath,colonequals,bm,nicefrac}
\usepackage{algorithm}
\usepackage{algorithmicx,comment}
\usepackage[noend]{algpseudocode}
\usepackage{tkz-graph}
\usepackage{url}

\newtheorem{theorem}{Theorem}[section]
\newtheorem{lemma}[theorem]{Lemma}

\newcommand{\R}[0]{{\mathbb{R}}}

\newcommand{\Ea}{{E_\mathsf a}}
\newcommand{\Eb}{{E_\mathsf b}}
\newcommand{\na}{{\nu_\mathsf a}}
\newcommand{\nb}{{\nu_\mathsf b}}

\newcommand{\Pad}{{\mathscr P_\mathsf{ad}}}
\newcommand{\ud}{{\bu_\mathrm d}}
\newcommand{\uo}{{\bu^\mathrm{obs}}}
\newcommand{\calC}{{\mathcal C}}
\newcommand{\Gr}{{\Gamma_\mathrm r}}
\newcommand{\Gl}{{\Gamma_\mathrm l}}
\newcommand{\Glr}{{\Gamma_\mathrm{lr}}}
\newcommand{\Gtb}{{\Gamma_\mathrm{tb}}}
\newcommand{\Slr}{{\Sigma_\mathrm{lr}}}
\newcommand{\Stb}{{\Sigma_\mathrm{tb}}}
\newcommand{\bu}{{\bm u}}
\newcommand{\bv}{{\bm v}}
\newcommand{\bVd}{{\bm V_\mathrm d}}
\newcommand{\bVdh}{{\bm V^h_\mathrm d}}
\newcommand{\bVh}{{\bm V^h}}
\newcommand{\bph}{{\bm\varphi}}
\newcommand{\bps}{{\bm\psi}}
\newcommand{\bx}{{\bm x}}
\newcommand{\bn}{{\bm n}}
\newcommand{\bt}{{\bm \tau}}
\newcommand{\yai}{{y_\mathsf a^i}}
\newcommand{\ybi}{{y_\mathsf b^i}}

\def\eop{{\ \vrule height 7pt width 7pt depth 0pt}}

\definecolor{DarkRed}{rgb}{0.5 0 0}
\definecolor{DarkBlue}{rgb}{0 0 0.5}
\definecolor{DarkGreen}{rgb}{0 0.5 0}
\definecolor{DarkOrange}{rgb}{0.83 0.33 0}
\definecolor{DarkMagenta}{rgb}{0.83 0 0.67}

\newcommand{\sfrei}{}

\usepackage{newtxtext}       %
\usepackage[varvw]{newtxmath}       

\makeindex             

\begin{document}

\title{A non-intrusive neural-network based BFGS algorithm for parameter estimation in non-stationary elasticity}
\author{Stefan Frei,
Jan Reichle
and Stefan Volkwein~\thanks{Department of Mathematics and Statistics, University of Konstanz, 78457 Konstanz, Germany, \{stefan.frei, jan.reichle, stefan.volkwein\}@uni-konstanz.de} }


%
\maketitle


\begin{abstract}
We present a non-intrusive gradient \sfrei{and a non-intrusive BFGS} algorithm for parameter estimation problems in non-stationary elasticity. To avoid multiple (and potentially expensive) solutions of the underlying partial differential equation (PDE), we approximate the PDE solver by a neural network within the \sfrei{algorithms}. The network is trained offline for a given set of parameters. The \sfrei{algorithms are} applied to an unsteady linear-elastic contact problem; their convergence and approximation properties are investigated numerically.
\end{abstract}



\section{Introduction}
\label{sec:1}

We are interested in the solution of parameter estimation problems subject to computationally expensive partial differential equations (PDE). Typically, the available algorithms to solve such problems require multiple solutions of the underlying PDE to compute state and sensitivities. To fix ideas, we study a model from elasticity describing the deformation of a metal plate under applied forces, where the material parameters $p\in {\cal P}^{\text{ad}}$ are unknown and have to be estimated. In industrial practice commercial software is frequently used to predict, for instance, the resulting deformations and stresses. In many cases, the software is not accessible and does not allow for a computation of sensitivities with respect to the parameters. 

In this work we apply a neural network to approximate the solution of a non-stationary linear-elastic contact problem. For an overview of data-based approaches to solve inverse problems, we refer, e.g., to~\cite{Arridgeetal2019}, but see also \cite{AEOV23,RPK19}. The approximation quality of neural networks is studied in~\cite{Grohsetal2023}, for instance, and the many references cited therein.

The paper is organized as follows: In Section~\ref{Section:2} the parameter identification problem is introduced as a PDE-constrained minimization problem. The finite element (FE) discretization and the neural network are briefly explained in Sections~\ref{Section:3} and \ref{sec.NN}, respectively. Section~\ref{Section:5} is devoted to the gradient- \sfrei{and the BFGS-} based optimization method, where the derivatives are computed by backpropagation. In Section~\ref{Section:6} a numerical example is presented. Finally, a conclusion is made in Section~\ref{Section:7}.


\section{Problem Formulation}
\label{Section:2}

\begin{figure}[t]
\centering
\includegraphics[width=0.8\textwidth]{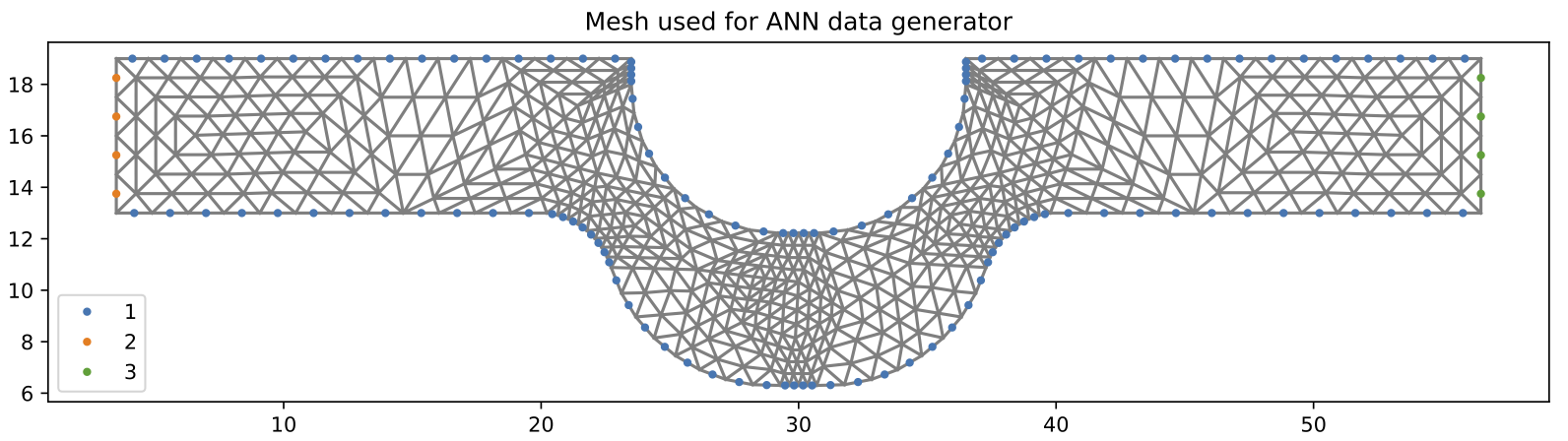}
 	\caption{Spatial domain $\Omega$ and its finite element triangulation $\mathscr T_h$.}
 	\label{fig:MeshNet}
\end{figure}

In the present work we consider a non-stationary linear-elastic contact problem, see, e.g., \cite{ChoulyHild2013}, in a time interval $[0,T]$ and in a two-dimensional bounded spatial domain $\Omega$. However, our approach can also be applied to other evolution equations. We assume that a plate lies in between a fixed rigid geometry in the form of a ridge, see Figure~\ref{fig:MeshNet}.

The plate is pulled from the right boundary ($\Gr$), while it is fixed on the left one ($\Gl$). On top and bottom ($\Gtb$) we apply frictionless contact conditions which allow the plate to deform only towards its interior. Then, we have $\partial\Omega = \Gtb \cup \Glr$, where $\Glr=\Gl\cup\Gr$. Moreover, we set
\begin{align*}
    \bm H=L^2(\Omega;\mathbb R^2)\quad\text{and}\quad\bVd=\big\{\bph\in H^1(\Omega;\mathbb R^2)\,\big|\,\bph=0\text{ on }\Glr\big\}.
\end{align*}
For $T>0$ let $Q=(0,T)\times\Omega$, $\Stb=(0,T)\times\Gtb$ and $\Slr=(0,T)\times\Glr$. The resulting deformation vector $\bu=(u_1,u_2)$ is governed by the initial-boundary value problem
\begin{subequations}
    \label{Eq:ModelProblem}
	\begin{align} 
		\rho\bu_{tt}-\nabla\cdot\sigma(\bu)&=\bm f&&\text{in }Q,\\
		\bu\cdot\bn \le 0,\hspace{1mm}\bn\cdot\big(\sigma(\bu)\bn\big) &\le 0,\hspace{1mm}(\bu\cdot \bn)\big((\bn\cdot(\sigma(\bu)\bn)\big) =0 &&\text{on }\Stb,\label{contact1}\\
		\bt\cdot\big(\sigma(\bu)\bn\big) &=0&&\text{on }\Stb,\\\label{contact2}
		\bu&=\ud&&\text{on }\Slr,\\
		\bu(0,\cdot)=\bu_t(0,\cdot)&=0&&\text{on }\Omega
	\end{align}
\end{subequations}
with an inhomogeneity $\bm f\in L^2(0,T;\bm H)$, the matrices $\sigma(\bu)=2\mu\epsilon(u)+\lambda\mathrm{tr}(\epsilon(\bu))I$ and $\epsilon(\bu)=\nicefrac{1}{2}\,(\nabla\bu+(\nabla\bu)^\top)$, where $I$ stands for the $2\times 2$ identity matrix. The (positive) Lam\'e parameters $\lambda$ and $\mu$ correspond to the elasticity modulus $E$ and to the Poisson ratio $\nu$ (which will be the parameters we wish to determine by means of experiments) by the relations
{$\displaystyle 
	\lambda=\frac{E\nu}{(1+\nu)(1-2\nu)}$} and {$\displaystyle \mu= \frac{E}{2(1+\nu)}.$}
The convex and compact set of admissible parameters is given as 
\begin{align*}
    \Pad\colonequals\big\{p=(E,\nu)\in\R^2\,\big|\,\Ea\le E\le\Eb,~\na\le\nu\le\nb\big\}\equalscolon\big[\Ea,\Eb\big]\times\big[\na,\nb\big],
\end{align*}
where $0<\Ea\le\Eb$ and $0<\na\le\nb< \nicefrac{1}{2}$. The density is fixed to $\rho = 2700\,\nicefrac{\mathrm{kg}}{\mathrm m^3}$ and the (spatially independent) Dirichlet conditions are given by 
\begin{align*}
\ud= 0 \text{ on }(0,T] \times \Gl,\quad\ud(t) =(18t^2 -12t^3, 0)\text{ on }(0,T] \times \Gr.
\end{align*} 
In \eqref{contact1} and \eqref{contact2} the vectors $\bn$ and $\bt$ denote the outer unit normal vector and the corresponding tangential vector, respectively.

To deal with the contact conditions~\eqref{contact1} we use Nitsche's method \cite{AlartCurnier91}. It can be shown~\cite{AlartCurnier91} that the three conditions~\eqref{contact1} are equivalent to the following equality for arbitrary $\gamma>0$:
\begin{equation}
    \label{eq:gammagleichung}
	\bn\cdot\big(\sigma(\bu)\bn\big)=-\frac{1}{\gamma}\left[\bu\cdot\bn-\gamma\bn\cdot\big(\sigma(\bu)\bn\big) \right]_+,
\end{equation}
where $[g]_+\colonequals\max\{g,0\}$ denotes the positive part of any function $g$. Condition \eqref{eq:gammagleichung} can be easily incorporated in a variational formulation; cf. \cite{ChoulyHild2013}, for instance. For that purpose we introduce the velocity vector $\bv=\bu_t$. Then the variational formulation of \eqref{Eq:ModelProblem} is given as
\begin{align}
    \label{Varform}
    \begin{aligned}
        \frac{\mathrm d}{\mathrm dt}\,{\langle\bv,\bph\rangle}_{\bm H}+\int_\Omega \sigma(\bu):\nabla \bph\, \mathrm d\bx\hspace{20mm}&\\
        +\frac{1}{\gamma}\,\Big\langle\left[\bu\cdot\bn-\gamma\bn\cdot\big(\sigma(\bu)\bn\big)\right]_+, \bph\cdot\bn\Big\rangle_{L^2(\Gtb)}&= 0\quad\forall\bph\in\bVd\text{ in }(0,T],\\
        \frac{\mathrm d}{\mathrm dt}\,{\langle \bu,\bps\rangle}_{\bm H}-{\langle\bv,\bps\rangle}_{\bm H} &= 0\quad\forall\bps\in\bm H\text{ in }(0,T],
    \end{aligned}
\end{align}
where the symbol ``$\,:\,$'' stands for the Frobenius product of matrices.

For known results on the well-posedness of dynamic contact problems, we refer to~\cite[Chapter 4]{Ecketal2005}. Throughout, we assume that the solution $\bu$ to \eqref{Varform} belongs to $C(\overline Q)$.

Next we formulate the parameter identification as an optimization problem. Suppose we are given five spatial points $\bx_i=(x_i,y_i)\in\Omega$ for $i=1,\ldots,5$. We set $t_1=\nicefrac{T}{2}$, $t_2=T$ and $\yai\colonequals\inf\{y\in\R\,|\,(x_1^i,y)\in\Omega\}$, $\ybi\colonequals\sup\{y\in\R\,|\,(x_1^i,y)\in\Omega\}$ for $i=1,\ldots,5$.
%
Assume that we have the following available measurements:
\begin{itemize}
    \item [1)] point values 
    of the deformation vector
    \begin{align*}
        C^1(\bu)\colonequals\big(\bu(\bx_1,t_1),\bu(\bx_2,t_1),\ldots,\bu(\bx_5,t_1),\bu(\bx_1,t_2),\ldots,\bu(\bx_5,t_2\big)\in\R^{2\times10};
    \end{align*}
    \item [2)] mean values
    \begin{align*}
        C^2(\bu)&\colonequals\textstyle\Big(\kappa_1\int_{y_\mathsf a^1}^{y_\mathsf b^1}\bu(x_1,y,t_1)\,\mathrm dy,\ldots,\kappa_5\int_{y_\mathsf a^5}^{y_\mathsf b^5}\bu(x_5,y,t_1)\,\mathrm dy,\\
        &\hspace{7mm}\textstyle\kappa_1\int_{y_\mathsf a^1}^{y_\mathsf b^1}\bu(x_1,y,t_2)\,\mathrm dy,\ldots,\kappa_5\int_{y_\mathsf a^5}^{y_\mathsf b^5}\bu(x_5,y,t_2)\,\mathrm dy\Big)\in\R^{2\times10}
    \end{align*}
    with {the factor} $\kappa_i=\nicefrac{1}{(\ybi-\yai)}$;
    \item [3)] mean values
    \begin{align*}
        C^3(\bu)&\colonequals\textstyle\Big(\kappa_1\int_{y_\mathsf a^1}^{y_\mathsf b^1}\sigma_{\text{vM}} (x_1,y,t_1)\,\mathrm dy,\ldots,\kappa_5\int_{y_\mathsf a^5}^{y_\mathsf b^5}\sigma_{\text{vM}} (x_5,y,t_1)\,\mathrm dy,\\
        &\hspace{7mm}\textstyle\kappa_1\int_{y_\mathsf a^1}^{y_\mathsf b^1}\sigma_{\text{vM}} (x_1,y,t_2)\,\mathrm dy,\ldots,\kappa_5\int_{y_\mathsf a^5}^{y_\mathsf b^5}\sigma_{\text{vM}} (x_5,y,t_2)\,\mathrm dy\Big)\in\R^{1\times10}
    \end{align*}
    of the van-Mises stress 
    \begin{align*}
        \sigma_{\text{vM}}\colonequals\big(\nicefrac{3}{2}\,\sigma^{\text{dev}}(\bu) : \sigma^{\text{dev}}(\bu)\big)^{1/2}, \text{ where } \sigma^{\text{dev}}(\bu)\colonequals\sigma(\bu) -\nicefrac{1}{3}\,\mathrm{tr}(\sigma(\bu)) I.
    \end{align*}
\end{itemize}
\noindent
These $m\colonequals50$ measurements define the nonlinear observation operator
\begin{align*}
    \bu\mapsto \calC(\bu)={\big(C^1_{1,\cdot}(\bu),C^1_{2,\cdot}(\bu),C^2_{1,\cdot}(\bu),C^2_{2,\cdot}(\bu),C^3(\bu)\big)^\top}\in\R^m.
\end{align*}
For given measurement vector $\uo\in\R^m$ we consider the minimization problem
\begin{align}
    \label{func1}
	\begin{aligned}
	    &\min J(\bu,(E,\nu))\colonequals\frac{1}{2}\,{\|{\mathcal C}(\bu)-\uo\|}^2_{\R^m}+\frac{\alpha}{2}\,\bigg(\Big(\frac{E}{E_0})^2+\big(\frac{\nu}{\nu_0}\Big)^2\bigg)\\
     &\hspace{0.5mm}\text{subject to }(\bu,(E,\nu))\text{ satisfies~\eqref{Varform} and }(E,\nu)\in\Pad,
	\end{aligned}
\end{align}
where the parameters $E, \nu$ are normalized by characteristic values $E_0, \nu_0$ of elasticity modulus and Poisson ratio, respectively.

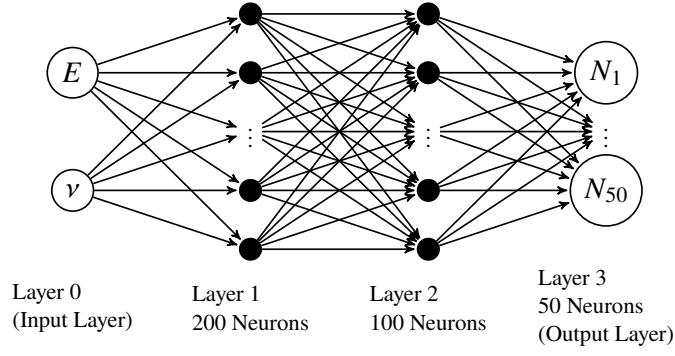
\begin{figure}[t]
\centering
 \begin{minipage}{0.6\textwidth}
	\scalebox{0.78}{\begin{tikzpicture}[->,>=stealth',shorten >=1pt, auto, node distance=2cm,thick,main node/.style={circle,draw,font=\Large\bfseries}]
        \node[main node] (1) at (0,0) {$E$};
        \node[main node] (2) at (0,-2) {$\nu$};
        \node[] (d1) at (3,-1) {$\vdots$};
        \node[] (d2) at (6,-1) {$\vdots$};
        \node[] (d3) at (9,-1) {$\vdots$};
        \node[main node, fill=black] (3) at (3,0) {};
        \node[main node, fill=black] (4) at (3,-2) {};
        \node[main node, fill=black] (3o) at (3,1) {};
        \node[main node, fill=black] (4o) at (3,-3) {};
        
        \node[main node, fill=black] (6) at (6,0) {};
        \node[main node, fill=black] (7) at (6,-2) {};
        \node[main node, fill=black] (5) at (6,1) {};
        \node[main node, fill=black] (8) at (6,-3) {};
        \node[main node] (9) at (9,0) {$N_1$};
        \node[main node] (10) at (9,-2) {$N_{50}$};
  
        \node [anchor=base east,align=left] (14) [below of=2] {Layer 0\\(Input Layer)};
        \node [anchor=base east,align=left] (15) [below of=4] {Layer 1\\$200$ Neurons};
        \node [anchor=base east,align=left] (16) [below of=7] {Layer 2\\$100$ Neurons};
        \node [anchor=base east,align=left] (17) [below of=10] {Layer 3\\$50$ Neurons\\(Output Layer)};
        \path
    	(1) 	edge node {} (3)
         	      edge node {} (4)
                edge node {} (d1)
                edge node {} (3o)
                edge node {} (4o)
    	(2) 	edge node {} (3)
         	      edge node {} (4)
                edge node {} (d1)
                edge node {} (3o)
                edge node {} (4o)
        (d1) 	edge node {} (5)
    		      edge node {} (6)
		          edge node {} (7)
		          edge node {} (8)
                edge node {} (d2)
    	(3) 	edge node {} (5)
    		      edge node {} (6)
		          edge node {} (7)
		          edge node {} (8)
                edge node {} (d2)
        (3o) 	edge node {} (5)
    		      edge node {} (6)
		          edge node {} (7)
		          edge node {} (8)
                edge node {} (d2)
    	(4) 	edge node {} (5)
    		      edge node {} (6)
		          edge node {} (7)
		          edge node {} (8)
                edge node {} (d2)
        (4o) 	edge node {} (5)
    		      edge node {} (6)
		          edge node {} (7)
		          edge node {} (8)
                edge node {} (d2)
        (d2) 	edge node {} (9)
    		      edge node {} (10)
		          edge node {} (d3)
	    (5) 	edge node {} (9)
         	      edge node {} (10)
                edge node {} (d3)
	    (6) 	edge node {} (9)
         	      edge node {} (10)
                edge node {} (d3)
	    (7) 	edge node {} (9)
         	      edge node {} (10)
                edge node {} (d3)
	    (8) 	edge node {} (9)
         	      edge node {} (10)
                edge node {} (d3)
	    (9) 	
	    (10) 	
	    ;      
	\end{tikzpicture}}
 \end{minipage}
 	\caption{Dense neural network with two hidden layers. \sfrei{We use the least-squares cost functional $\mathcal L(p^{(i)})=\frac{1}{2n}\sum_{i=1}^n{\|C^{(i)}-N(p^{(i)}) \|}_{\R^m}^2$ for given parameters $p^i=(E^i, \nu^i)$ and corresponding observations $C^{(i)}=C(u_i^h),\, i=1,\ldots,n$.}}
 	\label{fig:NN}
\end{figure}


\section{Finite Element (FE) Discretization}
\label{Section:3}

For discretization, we use the space of $P_1$ finite elements on a triangular mesh $\mathscr T_h$:
\begin{align*}
	V^h\colonequals\big\{v_h\in C(\overline{\Omega})\,\big|\, v_h|_T\in P_1(T)\, \forall T \in\mathscr T_h\big\},~\bVh=V^h\times V^h,~\bVdh\colonequals\bVh\cap\bVd.
\end{align*}
The FE mesh is illustrated in Figure~\ref{fig:MeshNet}, and the discrete formulation of \eqref{Varform} reads: \textit{Find }$\bu^h\in\ud+\bVdh$ and $\bv^h \in\bVh$, such that	$\bu^h(0)=\bv^h(0)=0$ in $\bm H$ and
\begin{align}
    \label{VarformFE}
    \begin{aligned}
        \frac{\mathrm d}{\mathrm dt}\,{\langle\bv^h,\bph^h\rangle}_{\bm H}+\int_\Omega \sigma(\bu^h):\nabla \bph^h\, \mathrm d\bx\hspace{15mm}&\\
        +\frac{1}{\gamma}\,\Big\langle\left[\bu^h\cdot\bn-\gamma\bn\cdot\big(\sigma(\bu^h)\bn\big)\right]_+, \bph^h\cdot\bn\Big\rangle_{L^2(\Gtb)}&= 0\quad\forall\bph^h\in\bVdh\text{ in }(0,T],\\
        \frac{\mathrm d}{\mathrm dt}\,{\langle \bu^h,\bps^h\rangle}_{\bm H}-{\langle\bv^h,\bps^h\rangle}_{\bm H} &= 0\quad\forall\bps\in\bVh\text{ in }(0,T].
    \end{aligned}
\end{align}
The Nitsche parameter is chosen as $\nicefrac{1}{\gamma}=\nicefrac{\gamma^0}{h}$ for a parameter $\gamma^0 > 0$, see e.g.,~\cite{ChoulyHild2013}. System~\eqref{VarformFE} is discretized in time using the backward Euler method.
					

\section{Neural Network Approximation}
\label{sec.NN}

To solve the {optimization} problem \eqref{func1} numerically, we have to compute solution pairs $(\bu^h,\bv^h)$ to \eqref{VarformFE} for many parameters $p=(E,\nu)\in\Pad$. For that reason we introduce a neural network with the goal to approximate the FE solution map $(\bu^h,\bv^h)=\mathcal S^h(p)$. We assume that \eqref{VarformFE} can be solved \textit{offline} for a certain number of parameter pairs $p^{(i)}=(E^{(i)},\nu^{(i)})$, $i=1,\ldots,n$, which is used to train the network.

We will use a fully connected dense network, which can be defined by the following map 
\begin{align*}
	N(W,B;p)=W^L(\sigma(W^{L-1}+\dots\sigma(W^1 p+b^1)\dots+b^{L-1}) +b^L,
\end{align*}
where $B = \{b^l\in\R^{n_l}\,|\,l\in \{1,\dots,L\}\}$ is a set of biases and $W=\{W^{l}\in\R^{n_l \times n_{l-1}}\,|\,l\in\{1,\dots,L\}\}$ a set of matrices, which are the parameters of the neural network. An illustration of a (small) network is shown in Figure~\ref{fig:NN}. The network receives a set $\{p^{(i)}\}_{i=1}^n\subset\Pad$ of parameters as input and shall deliver the observation $C^{(i)}=\mathcal C(\bu^h_i)\in\R^m$, $\bu^h_i=\mathcal S^h(p^{(i)})$, as an output. In the numerical examples, we will use the softplus activation function $\sigma(s)=\ln(e^s+1)$.

To train the network, we use a least-squares loss function, defined by
\begin{align*}
    \mathcal L(W,B)=\frac{1}{2n}\sum_{i=1}^n{\|C^{(i)}-N(W,B;p^{(i)}) \|}_{\R^m}^2.
\end{align*}
Both input and output data are normalized around the respective mean values and with respect to the variances before they are used to train the network, e.g., for the $j$-parameter $p_j$ by $\tilde p_j^{(i)}\colonequals\nicefrac{(p_j^{(i)}-\mathrm{mean}_j)}{\delta_j}$ with $\mathrm{mean}_j=\sum_{i=1}^np^{(i)}/n$ and variance $\delta_j$.
	            
\bigskip

\noindent \textbf{Training.}\,
We apply the neural network to the model problem \eqref{VarformFE}. To train the neural network, \eqref{VarformFE} has been solved for $100\,000$ equally distributed parameter pairs ($E,\nu$) in $[5\cdot 10^{10}, 10^{11}]\times [0.3,0.4]$ using $1000$ equally distributed values for $E$ and $100$ equally distributed values for $\nu$. 
For validation, we create another set of $19\,483$ data pairs corresponding to randomly created parameters $p\in\Pad$. We use the deep-learning API \textit{keras}, which is based on \textit{tensorflow}~\cite{chollet2015keras}. {A dense network with two hidden layers of size $200$ and $100$ is applied}. 
For the training, we use a combination of the stochastic gradient algorithm and the Adam algorithm, see~\cite{kingma2017adam}, for in total 500 epochs. More precisely, we execute both algorithms for 50 epochs each and use the result of the one, which performed better in terms of a reduction in the loss function evaluated \sfrei{with} validation data. We start with a learning rate of $10^{-3}$, which is increased gradually up to $5\cdot 10^{-2}$ in the last epochs. As a result the loss function decreased by a factor of roughly $10^4$ within 500 epochs. The total time required for training was approximately $107,6$ seconds.

\bigskip
		
\begin{figure}[t]
    \centering
    \includegraphics[width=5.6cm,height=4.8cm]{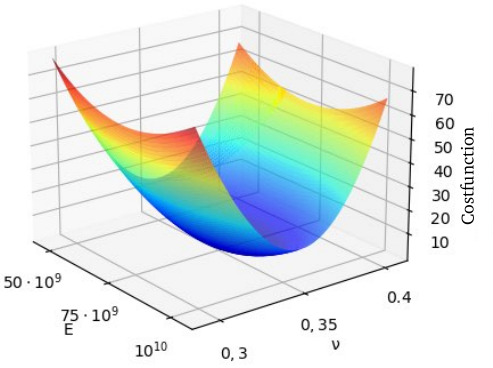} \hfill	
    \includegraphics[width=5.6cm,height=4.8cm]{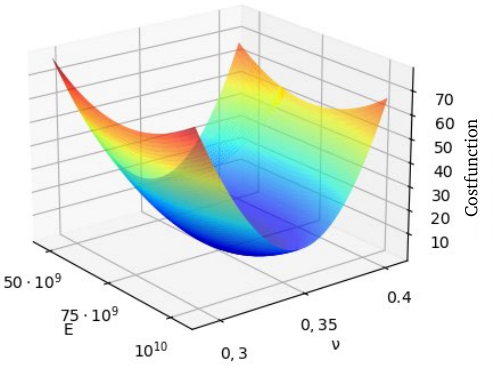}
	\caption{Comparison of the FE objective $\mathcal F^h(p_1)$ and its neural network approximation $\mathcal F^N(p_1)$ for $\uo=\mathcal C(\bu_1^h)$ and $\bu^h_1=\mathcal S^h(p_1)$ with $p_1$ defined in \eqref{punkt1}.
	\label{P1_cost_big}}
\end{figure}

\noindent\textbf{Approximation properties.}\,
In order to test the neural network approximation we compare the two functionals
\begin{align}
    \label{ges_min_funk}
    \begin{aligned}
        \mathcal F^h(E,\nu)&\colonequals\frac{1}{2}\,{\|{\mathcal C}(\bu^h)-\uo\|}^2_{\R^m}&&\text{for }\bu^h=\mathcal S^h(E,\nu),\\
        \mathcal F^N(E,\nu)&\colonequals\frac{1}{2}\,{\|{\mathcal C}(\bu^N)-\uo\|}^2_{\R^m}&&\text{for }\bu^N=N(W,B;(E,\nu)).
    \end{aligned}
\end{align}
We assume that the parameters we would like to find in our parameter estimation problem are given by
\begin{align}
    \label{punkt1}
	p_1\colonequals(E_1,\nu_1)\quad\text{with }E_1=7.513\cdot 10^9,~\nu_1=0.3547
\end{align}
and define the observed data by means of the FE solution $\bu^h_1=\mathcal S^h(p_1)$ of \eqref{VarformFE} as 
$\uo=\mathcal C(\bu_1^h)$.
In Figure~\ref{P1_cost_big}, we compare $\mathcal F^h$ and its neural network approximation $\mathcal F^N$ on $\Pad$. We observe that both graphs nearly coincide. Moreover, $\mathcal F^N$ seems to be convex despite the non-linear activation function $\sigma$ of the network, which is beneficial for numerical optimization. The minimum is attained very close to $p_1$, see also the numerical results for the parameter estimation problem given below.

\section{Neural-Network Based Gradient \sfrei{and BFGS} Algorithms}
\label{Section:5}

Gradient-type algorithms typically require frequent solutions of the PDE \eqref{VarformFE} and its sensitivities. In order to avoid these, we approximate the operator $p \mapsto \mathcal C(\mathcal S^h(p))$) by means of a neural network, i.e., we seek the minimizer of the function $\mathcal F^N(p)$ defined in \eqref{ges_min_funk}. Therefore, we need to compute the derivative of the network output $N(W,B;p)$ with respect to $p\in\Pad$. This can be computed by means of a backpropagation: 

\begin{lemma}\label{lem.grad}
    Let $N(W,B;p)$ be the dense network defined in Section~{\em\ref{Section:3}}. Moreover, let $z^{l}=W^{l}a^{l-1}+b^{l}$ and $a^{l}=\sigma_l(z^{l})$ be
	the output of the $l$-the layer of the network. It holds for the derivatives $\rho^l=(\rho^{l}_j)_{j=1}^{n_l}$, $\rho^{l}_j=\partial_{a^{l}_j}\mathcal F^N$:
    \begin{align*}
        \rho^L = a^L - y \quad \text{ and }\quad 
		\rho^l =  (W^{l+1})^\top\big(\rho^{l+1} \circ \sigma_{l+1}'(z^{l+1})\big)\, \text{ for } l=0,\dots , L-1,
    \end{align*}
	where $(u\circ v)_j=u_jv_j$. In particular, we obtain the derivative $\rho_j^{0}\colonequals\partial_{p_j}\mathcal F^N$ of $\mathcal F^N$ with respect to the network input $p=(p_1, p_2)$.
\end{lemma}

\textit{Proof}.
    The proof for $l\geq 1$ is standard and can be found in any textbook on neural networks, see, e.g., \cite{higham2018deep}. The extension to $l=0$ is completely straight-forward.\hfill\eop
    \bigskip

We use Lemma~\ref{lem.grad} to \sfrei{define} a gradient algorithm to find the minimum of $\mathcal F^N(p)$. The resulting algorithm is given as Algorithm~\ref{Algo_Grad_meth1}.
For the line search, we use an Armijo-Goldstein line search in the numerical results given below.
\begin{algorithm}[ht]
		\caption{(Neural-network based gradient algorithm)}\label{Algo_Grad_meth1}
		\begin{algorithmic}[1]
			\Procedure{Gradient algorithm}{network $N$, $p^0$}
			\State $p \gets p^0$
			\While {Termination criterion is not satisfied}
				\State $\Delta p \gets$ calculate $\nabla_p \mathcal F^N(p)$ \hfill \% using Lemma~\ref{lem.grad}
				\State $\eta \gets$ Line\_Search$(\mathcal F^N(\cdot), p, \Delta p)$ 
				\State $p \gets p-\eta \Delta p$ 
			\EndWhile
			\Return $p$
		\EndProcedure
\end{algorithmic}
\end{algorithm}

\sfrei{The gradient algorithm can be improved by using (approximate) second-order derivative information. In Algorithm~\ref{Algo_BFGS_meth}, we present a neural-network based BFGS algorithm.}

\begin{algorithm}[ht]
		\caption{\sfrei{(Neural-network based BFGS algorithm)}}\label{Algo_BFGS_meth}
		\begin{algorithmic}[1]
			\Procedure{BFGS algorithm}{network $N$, $p^0$}
			\sfrei{\State $p \gets p^0$, \,$B \gets \mathbb{I}_2$
			\While {Termination criterion is not satisfied}
            \State Compute the direction $d = -B^{-1} \nabla_p \mathcal F^N(p)$
            \State $\eta \gets$ Line\_Search$(\mathcal F^N(\cdot), p, d)$
            \State $p \gets p-\eta d$
            \State Compute $s = -\eta d$ and $y = \nabla_p \mathcal F^N(p)- \nabla_p \mathcal F^N(p-s)$
            \If {$y^{\top} s \ge \varepsilon$}
                \State $B \gets B +\frac{y y^{\top}}{y^{\top} s} - \frac{Bs(Bs)^{\top}}{s^{\top} Bs}$
            \EndIf
			\EndWhile
			\Return $p$}
		\EndProcedure
\end{algorithmic}
\end{algorithm}

%


\section{Numerical Examples}
\label{Section:6}
\vspace*{-0.1cm}

We use again the example introduced in Section~\ref{sec.NN} with $\uo=\mathcal C(\bu^h_1)$ and $\bu_1^h=\mathcal S^h(p_1)$ for $p_1$ specified in \eqref{punkt1}. The neural network $N(W,B;\cdot)$ is trained \textit{offline} as described in Section~\ref{sec.NN}. An evaluation of the neural network $N(W,B;p)$ for a new parameter pair $p$ requires then only $29.8$ $\mu s$ in average, while our FE solver needs $243$ seconds to solve \eqref{VarformFE}.
Due to the convexity of the functional $\mathcal F^N$ (see Fig.~\ref{P1_cost_big}) a regularization is not required, i.e., we set $\alpha=0$ in~\eqref{func1}. \\

\begin{figure}[t]
	\centering
	\includegraphics[width=5.4cm,height=4.9cm]{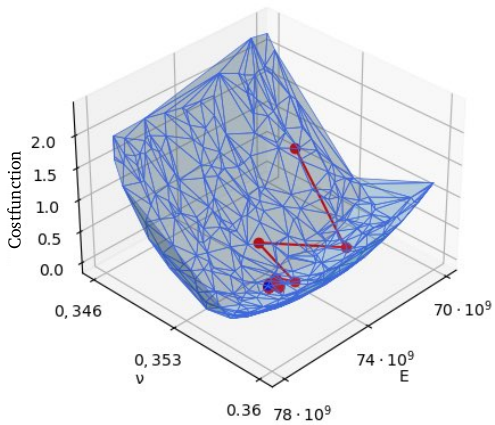}	\hfil		\includegraphics[width=5.4cm,height=4.9cm]{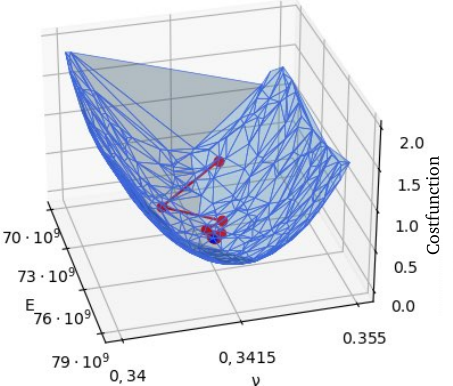}
	\caption{Convergence behavior of the neural-network based gradient algorithm applied to problem~\eqref{ges_min_funk} \sfrei{with exact solutions $p_1$~\eqref{punkt1} (left sketch) and $p_2$ (right sketch), respectively}. The correct solutions are marked with blue dots, the iterates in red.}
	\label{P1}
\end{figure}

\begin{figure}[t]
	\centering
	\includegraphics[width=4.8cm,height=4cm]{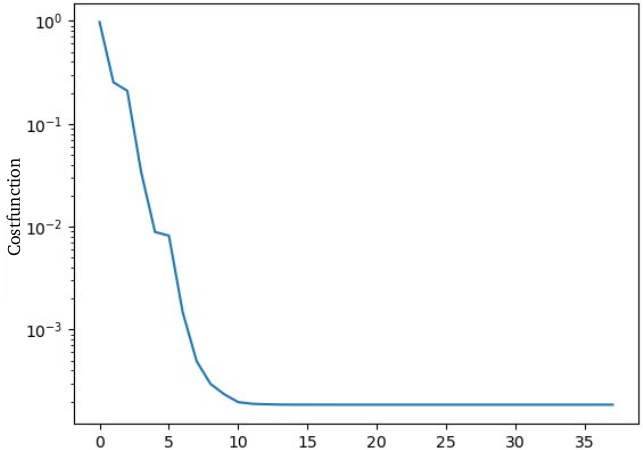}\hfil
	\includegraphics[width=4.8cm,height=4cm]{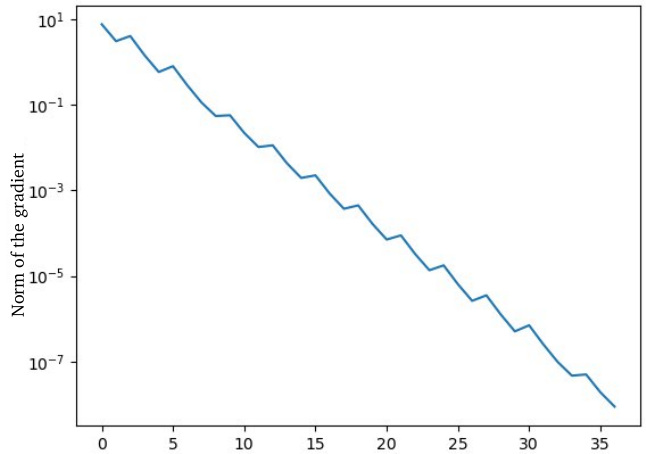}
	\caption{\textit{Left}: Objective function ${\cal F}^N$, \textit{right}: Norm of the gradient $\|\nabla_p {\cal F}^N\|_{\R^2}$ \sfrei{for the gradient algorithm applied to} the first example corresponding to parameters~\eqref{punkt1}.}
	\label{P1_cost_direct}
\end{figure}

\sfrei{\textbf{Gradient algorithm}}\,
\sfrei{First, we investigate the neural-network based gradient algorithm given in Algorithm~\ref{Algo_Grad_meth1}.
The algorithm is stopped,} when the initial value of the norm of the gradient $\|\nabla_p {\cal F}^N(p)\|_{\R^m}$ is reduced by a factor of $10^{-8}$. The {behavior} of the gradient algorithm is visualized in Figure~\ref{P1}, left sketch. We observe that the algorithm converges quickly towards the minimizer. This is also confirmed in Figure~\ref{P1_cost_direct}, where the objective function $\mathcal F^N$ and its gradient are visualized within the iteration. \sfrei{36 iterations are required until the stopping criterion is fulfilled.}
The resulting parameters $p$ are (after re-normalisation) $\bar{E}_1\approx 7.51385\cdot 10^{10}$, $\bar{\nu}_1\approx 0.354768$, which yields a relative error of approximately $1.1 \cdot 10^{-3}$. This seems to be the error resulting from the network approximation.
		
\sfrei{To confirm the results}, we apply the algorithm in the same way to a second point
%
	$p_2 = (E_2, \nu_2)$, where $E_2 = 7.45\cdot 10^{10}$ and $\nu_2 =0.3481$.
%
The convergence behavior is shown in Figure~\ref{P1}, right sketch. The resulting parameters are $\bar{E}_2 \approx 7.45071\cdot 10^{10}$,\, $\bar{\nu}_2\approx  0.348164$, which gives a relative error of around $3.8\cdot 10^{-3}$.\\

\begin{figure}[t]
	\centering
	\includegraphics[width=5.4cm,height=4.9cm]{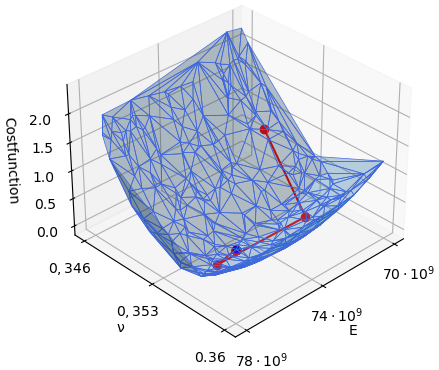}	\hfil		\includegraphics[width=5.4cm,height=4.9cm]{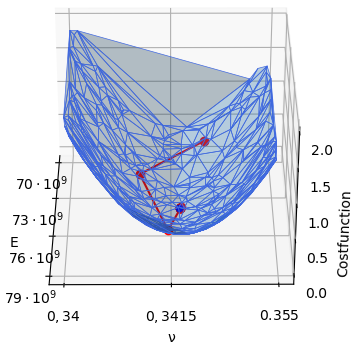}
	\caption{\sfrei{Convergence behavior of the neural-network based BFGS algorithm applied to problem~\eqref{ges_min_funk} with exact solutions $p_1$~\eqref{punkt1} (left sketch) and $p_2$ (right sketch), respectively.}}
	\label{P1_BFGS}
\end{figure}

\begin{figure}[t]
	\centering
	\includegraphics[width=4.8cm,height=4cm]{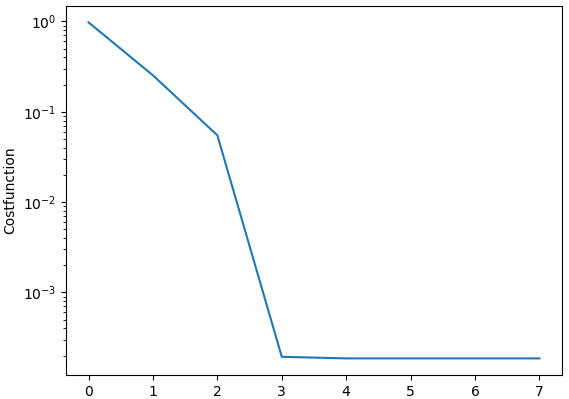}\hfil
	\includegraphics[width=4.8cm,height=4cm]{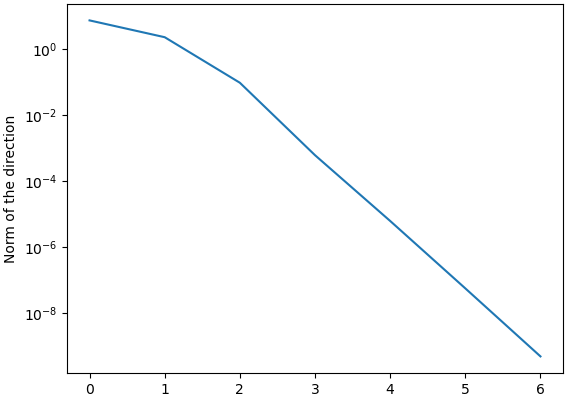}
	\caption{\sfrei{\textit{Left}: Objective function ${\cal F}^N$, \textit{right}: Norm of the gradient $\|\nabla_p {\cal F}^N\|_{\R^2}$ for the BFGS algorithm applied to the first example corresponding to parameters ~\eqref{punkt1}.}}
	\label{P1_cost_BFGS}
\end{figure}

\sfrei{\textbf{BFGS algorithm}}\,
\sfrei{Secondly, we apply the BFGS algorithm defined in Algorithm~\ref{Algo_BFGS_meth} to the two examples introduced above. The convergence behavior of the first few iterates is visualized in Figure~\ref{P1_BFGS}. We observe that the iterates converge in a more direct way towards the minimizer and most of the oscillations visible in Figure~\ref{P1} for the gradient algorithm are avoided. In Figure~\ref{P1_cost_BFGS}, we show the behavior of the objective function and its gradient over the iterates for the first example. The stopping criterion is already fulfilled after 7 iterations (compared to 36 iterations of the gradient algorithm), while the additional effort in each iteration is negligible. The minima found by the two algorithms are identical in the leading 8 digits.}

\vspace*{-0.3cm}		

\section{Conclusion}
\label{Section:7}
\vspace*{-0.2cm}

We have presented a non-intrusive gradient \sfrei{and a BFGS} algorithm based on a neural-network approximation of the underlying PDE. First results for a dynamic contact problem look very promising. \sfrei{The BFGS variant converges significantly faster compared to the gradient algorithm.} Despite the non-linear activation functions involved in the network approximation, the approximated cost function ${\cal F}^N$ seems to conserve the convexity of the function ${\cal F}^h$. We think that the algorithm can be particularly interesting in cases, when the PDE solvers are computationally very expensive or not accessible at all for the computation of sensitivities. 

\medskip

\bibliographystyle{plain}

\end{document}